\documentclass[11pt,leqno]{amsart}
\usepackage{indentfirst, amssymb,amsmath,amsthm}
\usepackage{hyperref}
\usepackage{csquotes}
\usepackage{graphicx}
\usepackage{epstopdf}
\usepackage{multicol}
\usepackage[margin=.75in]{geometry}
\newtheorem{theo}{Theorem}
\newtheorem{cor}{Corollary}
\newtheorem{rem}{Remark}
\newcommand{\beas}{\begin{eqnarray*}}
\newcommand{\eeas}{\end{eqnarray*}}
\begin{document}
\title[A visual tour via Definite Integration]{A visual tour via the Definite Integration $\int_{a}^{b}\frac{1}{x}dx$}
\date{}
\author[B. Chakraborty]{Bikash Chakraborty}
\date{}
%\footnotetext{MSC 2020: Primary: 11A99, 00A05, 11B33, 00A66.}
\address{Department of Mathematics, Ramakrishna Mission Vivekananda
Centenary College, Rahara, West Bengal 700 118, India.}
\email{bikashchakraborty.math@yahoo.com, bikash@rkmvccrahara.org}
\maketitle
\begin{abstract} Geometrically, $\int_{a}^{b}\frac{1}{x}dx$ means the area under the curve $\frac{1}{x}$ from $a$ to $b$, where $0<a<b$, and this area gives a positive number. Using this area argument, in this expository note, we present some visual representations of some classical results. For examples, we demonstrate an area argument on a generalization of Euler's limit $\left(\lim\limits_{n\to\infty}\left(\frac{(n+1)}{n}\right)^{n}=e\right)$. Also, in this note, we provide an area argument of the inequality $b^a < a^b$, where $e \leq a< b$, as well as we provide a visual representation of an infinite geometric progression. Moreover, we prove that the Euler's constant $\gamma\in [\frac{1}{2}, 1)$ and the value of $e$ is near to $2.7$.\par
Some parts of this expository article has been accepted for publication in Resonance – Journal of Science Education, The Mathematical Gazette, and International Journal of Mathematical Education in Science and Technology.
\end{abstract}
\section{Introduction}
It is well known that the function $\phi:(0,+\infty)\to \mathbb{R}$, defined by $\phi(x)=\frac{1}{x}$, is a monotone decreasing and continuous. Thus $\phi(x)$ is Riemann integrable on $[a,b]$ where $0<a<b$. Geometrically, $\int_{a}^{b}\frac{1}{x}dx$ means the area under the curve $y=\frac{1}{x}$ from $a$ to $b$. Moreover, it is useful to observe that the funtion $f(t)=\int_{1}^{t} \frac{1}{x}dx$ is strictly increasing for $t\geq 1$.\par
\begin{figure}[h]
\includegraphics[scale=.6]{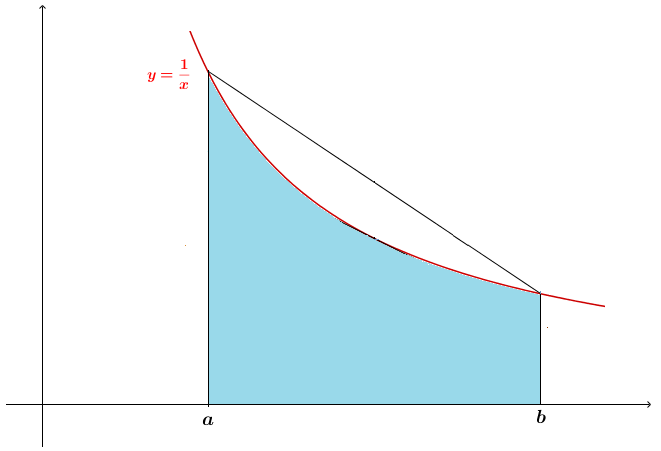}
\caption{$\ln b-\ln a< \frac{1}{2}\cdot(\frac{1}{a}+\frac{1}{b})\cdot(b-a).$}
\centering
\end{figure}
 Let $a$ and $b$ be two positive real numbers. Then the fact $\frac{1}{x}+\frac{x}{ab}\leq \frac{1}{a}+\frac{1}{b}$ for $a\leq x \leq b$ (as $(x-b)(x-a)\leq 0$)is equivalent to saying that the line $y=\frac{1}{a}+\frac{1}{b}-\frac{x}{ab}$ lies above the curve $y=\frac{1}{x}$ for $a\leq x \leq b$. Thus, Figure 1 shows that the area under the curve $y =\frac{1}{x}$ from $a$ to $b$ is less than the area of the trapeziums covering it, i.e., $$\ln b-\ln a< \frac{1}{2}\cdot(\frac{1}{a}+\frac{1}{b})\cdot(b-a).$$
 \begin{figure}[h]
\includegraphics[scale=.5]{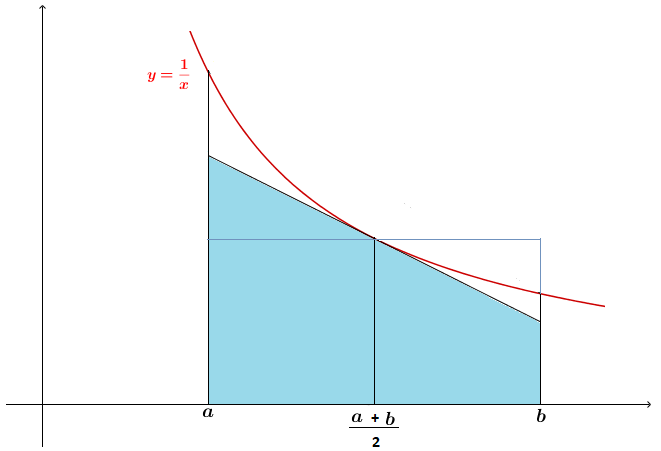}
\caption{$\ln b-\ln a>\frac{2}{a+b}\cdot(b-a).$}
\centering
\end{figure}
Again, the fact $\frac{1}{x}+\frac{4x}{(a+b)^{2}}\geq \frac{4}{a+b}$ for $x>0$ (follows from AM-GM inequality) is equivalent to saying that the curve $y=\frac{1}{x}$ lies above its tangent line $y=\frac{4}{a+b}-\frac{4x}{(a+b)^{2}}$ at the point $(\frac{a+b}{2},\frac{2}{a+b})$.\par Thus, figure 2 gives the visualization that the area under the
curve $y =\frac{1}{x}$ from $a$ to $b$ is greater than the area of the trapezium below it, i.e., $$\ln b-\ln a>\frac{2}{a+b}\cdot(b-a).$$
%These understandings are the main tools for our visual tour.
%%%%%%%%%%%%%%%%%%%%%%%%%%%%%%%%%%%%%%%%%%%%%%%%%%%%%%%%%%%%%%%%%%%%%%%%%%%%%%%%%%%%%%%%%%%%%%%%%%%%%%%%%%%%%%%%%%%%%%%%%%%%%%%%%%%%%%%%%%%%%%%%%%%%%%%%%%%%
\section{Tour-1}
In a recent note of the American Mathematical Monthly, R. Farhadian (\cite{i}) made a beautiful generalization of Euler's limit $\left(\lim\limits_{n\to\infty}\left(\frac{(n+1)}{n}\right)^{n}=e\right)$ as follows:
\begin{theo} (\cite{i})
Let $A_{n}$ be a strictly increasing sequence of positive numbers satisfying the asymptotic formula $A_{(n+1)}\sim A_{n}$, and let $d_{n}=A_{(n+1)}-A_{n}$. Then
$$\lim\limits_{n\to\infty}\left(\frac{A_{(n+1)}}{A_{n}}\right)^\frac{A_{n}}{d_{n}}=e.$$
\end{theo}
Now, we will provide a second proof of it, which is purely pictorial.
\begin{figure}[h]
\includegraphics[scale=.4]{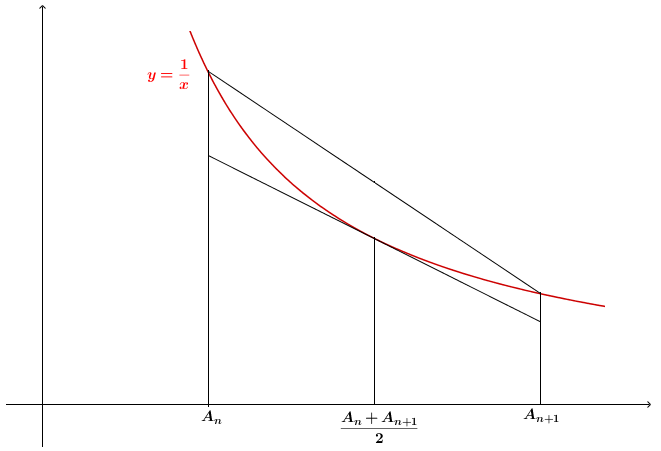}
\caption{}
\centering
\end{figure}
From  Figure 3, it is clear that \beas \frac{2}{A_{n}+A_{n+1}}\cdot(A_{n+1}-A_{n})<\ln(A_{n+1})-\ln(A_{n})<\frac{1}{2}\cdot\left(\frac{1}{A_{n}}+\frac{1}{A_{n+1}}\right)\cdot(A_{n+1}-A_{n}),\eeas
i.e., $$\frac{2}{1+\frac{A_{n+1}}{A_{n}}}<\ln\left(\frac{A_{n+1}}{A_{n}}\right)^{\frac{A_n}{d_n}}<\frac{1}{2}\cdot\left(1+\frac{A_{n}}{A_{n+1}}\right).$$
Since $A_{(n+1)}\sim A_{n}$, thus $\lim\limits_{n\to\infty}\left(\frac{A_{(n+1)}}{A_{n}}\right)^\frac{A_{n}}{d_{n}}=e.$
\begin{rem}
It is well-known that if  $a_{n}$  is a sequence of positive numbers satisfying $\lim\limits_{n\to +\infty} a_{n}=0$, then $$\lim\limits_{n\to +\infty} \left(1+a_{n}\right)^{\frac{1}{a_n}}=e$$.\\
Here, we will provide a visual proof of it.
\begin{figure}[h]
\includegraphics[scale=.4]{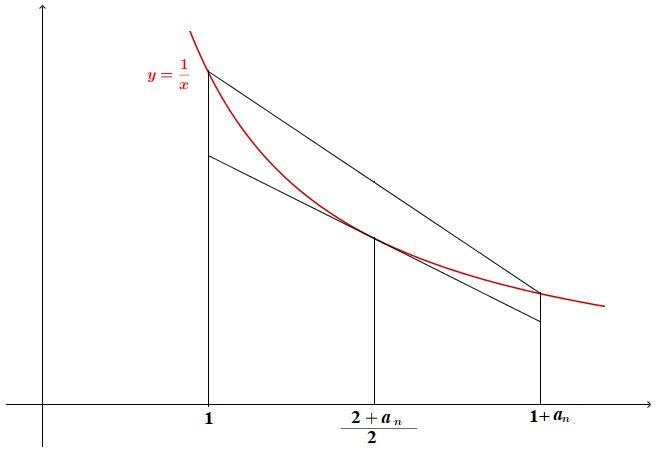}
\caption{}
\centering
\end{figure}
From the  figure, it is clear that \beas \frac{2}{2+a_{n}}\cdot a_{n}<\ln(1+a_{n})-\ln 1<\frac{1}{2}\cdot\left(1+\frac{1}{1+a_{n}}\right)\cdot a_{n},\eeas
i.e., $$\lim\limits_{n\to +\infty} \left(1+a_{n}\right)^{\frac{1}{a_n}}=e$$.
\end{rem}
%%%%%%%%%%%%%%%%%%%%%%%%%%%%%%%%%%%%%%%%%%%%%%%%%%%%%%%%%%%%%%%%%%%%%%%%%%%%%%%%%%%%%%%%%%%%%%%%%%%%%%%%%%%%%%%%%%%%%%%%%%%%%%%%%%%%%%%%%%%%%%%%%%%%%%%%%%%%
\section{Tour-2}
Next, we provide a pictorial description of a geometric series $$1+\frac{1}{r}+\frac{1}{r^2}+\frac{1}{r^3}+\ldots=\frac{r}{r-1}$$ when $r>1$.
\begin{figure}[h]
\includegraphics[scale=.2]{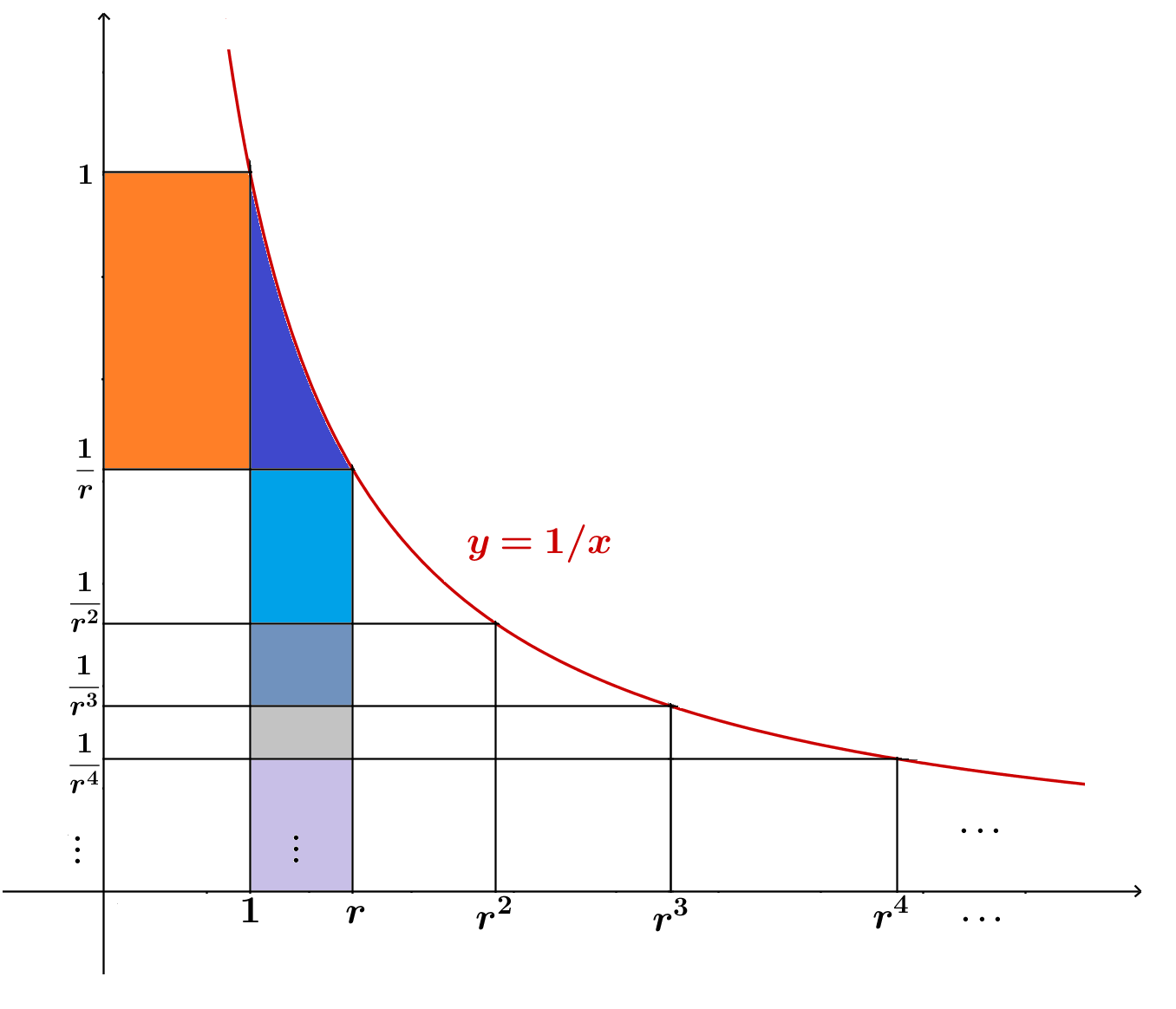}
\caption{$\int_{1}^{r}  \frac{dt}{t}=\int_{\frac{1}{r}}^{1}  \frac{dt}{t}$}
\centering
\end{figure}
Since $\int_{1}^{r}  \frac{dx}{x}=\int_{\frac{1}{r}}^{1}  \frac{dy}{y}$, thus the area covered by the rectangle $\{(0,1);(1,1);(1,\frac{1}{r});(0,\frac{1}{r})\}$ is same as the area covered by the rectangle $\{(1,0);(r,0);(r,\frac{1}{r});(1,\frac{1}{r})\}$. Thus
\beas \left(1-\frac{1}{r}\right)\cdot1&=&\left(\frac{1}{r}-\frac{1}{r^2}\right)\cdot(r-1)+\left(\frac{1}{r^2}-\frac{1}{r^3}\right)\cdot(r-1)+\ldots,\\
&=&\frac{(r-1)^2}{r^2}\cdot\left(1+\frac{1}{r}+\frac{1}{r^2}+\ldots\right),\eeas
i.e.,
$$\left(1+\frac{1}{r}+\frac{1}{r^2}+\ldots\right)=\frac{r}{r-1},$$ which gives the required equality.
%%%%%%%%%%%%%%%%%%%%%%%%%%%%%%%%%%%%%%%%%%%%%%%%%%%%%%%%%%%%%%%%%%%%%%%%%%%%%%%%%%%%%%%%%%%%%%%%%%%%%%%%%%%%%%%%%%%%%%%%%%%%%%%%%%%%%%%%%%%%%%%%%%%%%%%%%%%%%%%
\section{Tour-3}
The two constants $e$ and $\pi$ have encouraged many visual proofs of the inequality ${\pi}^e < e^{\pi}$. In a recent Mathematical Intelligencer note (\cite{C}), the author provided a visual proof of the inequality $\pi^{e}< e^{\pi}$. However, their visual proof can be used to show the more general inequality $b^a < a^b$, where $e \leq a< b$.
\medbreak
\textbf{Visual Proof-1}
\\
Since $\ln a\geq  1$, thus $\frac{1}{x\ln a}\leq \frac{1}{x}$ for $x>0$. Thus the Figure 6 shows that the area under the curve $y = \frac{1}{x\ln a}$ from $a$ to $b$ is less than  the area of the rectangle PQRS, i.e.,
\begin{figure}[h]
\includegraphics[scale=.6]{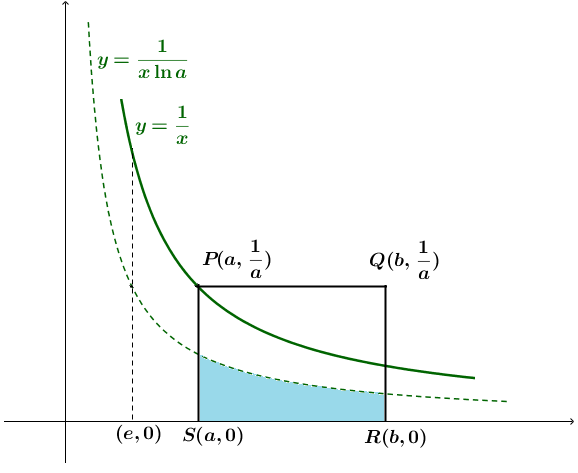}
\caption{$\int\limits_{a\ln a}^{b\ln a}\frac{dx}{x} < \frac{1}{a}\cdot(b-a)\ln a.$}
\centering
\end{figure}
\beas   \frac{\ln b }{\ln a}-1=\int_{a}^{b}\frac{dx}{x\ln a} < \frac{1}{a}(b-a)=\frac{b}{a}-1, \eeas
 % i.e.,   \beas\label{11}   a\ln b -a\ln a < (b-a)\ln a. \eeas
%Thus   \beas a\ln b<b\ln a, \text{~~i.e.,~~} b^{a}< a^{b}.\eeas
i.e., \beas b^{a}< a^{b}.\eeas
%%%%%%%%%%%%%%%%%%%%%%%%%%%%%%%%%%%%%%%%%%%%%%%%%%%%%%%%%%%%%%%%%%%%%%%%%%%%%%%%%%%%%%%%%%%%%%%%%%%%
\newpage
\textbf{Visual Proof-2}
\\
Also, Figure 7 shows that the area under the curve $y = \frac{1}{x}$ from $a\ln a$ to $b\ln a$ is less than  the area of the rectangle covering it. Since $e\leq a$, so $1\leq \ln a$, i.e., $a\leq a\ln a<b\ln a$.
\begin{figure}[h]
\includegraphics[scale=.6]{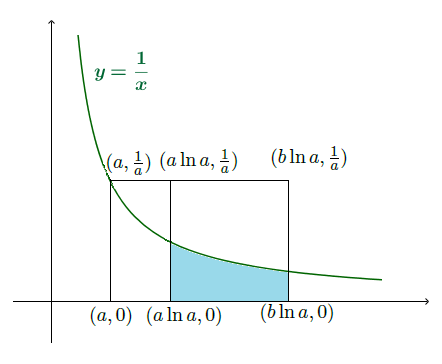}
\caption{$\int\limits_{a\ln a}^{b\ln a}\frac{dx}{x} < \frac{1}{a}\cdot(b-a)\ln a.$}
\centering
\end{figure}
\beas  \ln b -\ln a < \ln a\cdot\left(\frac{b}{a}-1\right),\eeas
i.e.,
 \beas   \frac{\ln b}{\ln a}-1 < \frac{b}{a}-1. \eeas
Thus \beas a\ln b<b\ln a,  \text{~~i.e.,~~} b^{a}< a^{b}.\eeas
\begin{cor} (\cite{C2})
If we take $a=e$, then $(a,0)$ and $(a\ln a, 0)$ will be coincided with $(e, 0)$. Also, $(a,\frac{1}{a})$ and $(a\ln a, \frac{1}{a})$ will be coincided with $(e, \frac{1}{e})$. Thus the figure 7 becomes :
\begin{figure}[h]
\includegraphics[scale=.6]{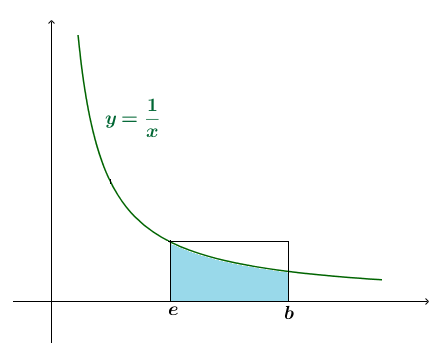}
\caption{$\int_{e}^{b}\frac{dx}{x} < \frac{1}{e}(b-e).$}
\centering
\end{figure}
\end{cor}
Thus we get $\ln b-1<\frac{b}{e}-1$, i.e., $b ^{e}<e^{b }$.
\begin{cor}
By taking $a=e$ and $b=\pi$, we get $\pi ^{e}<e^{\pi }$ (\cite{C}).
\end{cor}
%%%%%%%%%%%%%%%%%%%%%%%%%%%%%%%%%%%%%%%%%%%%%%%%%%%%%%%%%%%%%%%%%%%%%%%%%%%%%%%%%%%%%%%%%%%%%%%%%%%%%%%%%%%%%%%%%%%%%%%%%%%%%%%%%%%%%%%%%%%%%%%%%%%%%%%%%%
\section{Tour-4}
Considering the definition of the number $e$ by the equation $$1=\int_{1}^{e}\frac{1}{x} dx,$$ we are explaining that why the value of $e$ is near to $2.7$. Basically, we will show that $2.7<e<2.75$.
\begin{figure}[h]
\includegraphics[scale=.5]{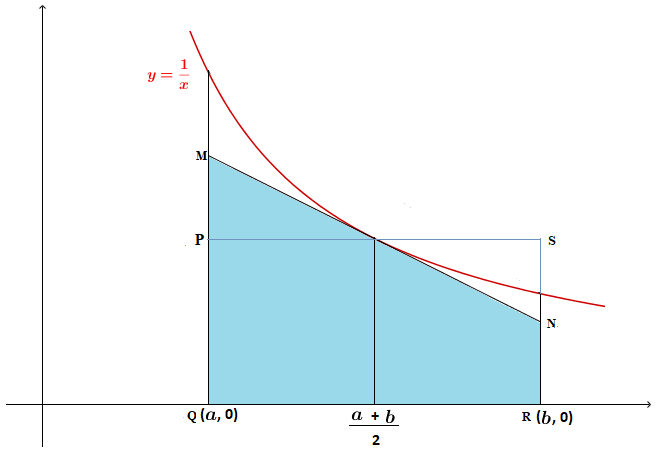}
\caption{$\int_{a}^{b} \frac{1}{x}dx > \frac{2(b-a)}{b+a}.$}
\centering
\end{figure}
Applying the lower bound of the integral $\int_{a}^{b} \frac{1}{x}dx$ (see,  Figure 9), we have \beas \int_{4}^{11} \frac{1}{x}dx &=& \int_{4}^{6} \frac{1}{x}dx +\int_{6}^{9} \frac{1}{x}dx +\int_{9}^{11} \frac{1}{x}dx\\
&>&\frac{2}{5}+\frac{2}{5}+\frac{1}{5}\\
&=&1, \eeas
i.e., $$ \int_{1}^{\frac{11}{4}} \frac{1}{x}dx > \int_{1}^{e}\frac{1}{x} dx~~~~~\Rightarrow~~e<2.75.$$
Again, applying the upper bound of the integral $\int_{a}^{b} \frac{1}{x}dx$ (see,  Figure 10), we have
\begin{figure}[h]
\includegraphics[scale=.5]{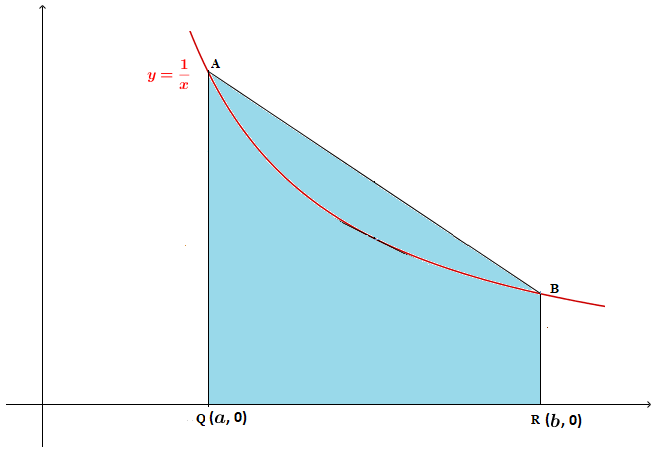}
\caption{$\int_{a}^{b} \frac{1}{x}dx < \frac{1}{2}\cdot(\frac{1}{a}+\frac{1}{b})\cdot(b-a).$}
\centering
\end{figure}

\beas &&\int_{10}^{27} \frac{1}{x}dx\\ &=& \int_{10}^{12} \frac{1}{x}dx +\int_{12}^{15} \frac{1}{x}dx +\int_{15}^{18} \frac{1}{x}dx+\int_{18}^{21}\frac{1}{x}dx+\int_{21}^{24} \frac{1}{x}dx +\int_{24}^{27} \frac{1}{x}dx\\
&<&\frac{1}{2}\cdot(\frac{2}{10}+\frac{2}{12})+\frac{1}{2}\cdot(\frac{3}{12}+\frac{3}{15})+\frac{1}{2}\cdot(\frac{3}{15}+\frac{3}{18})+
+\frac{1}{2}\cdot(\frac{3}{18}+\frac{3}{21})\\&&+\frac{1}{2}\cdot(\frac{3}{21}+\frac{3}{24})+\frac{1}{2}\cdot(\frac{3}{24}+\frac{3}{27})\\
&<& 1, \eeas
i.e., $$\int_{1}^{2.7} \frac{1}{x}dx <\int_{1}^{e} \frac{1}{x}dx~~~~~~\Rightarrow~~ e>2.70.$$
%%%%%%%%%%%%%%%%%%%%%%%%%%%%%%%%%%%%%%%%%%%%%%%%%%%%%%%%%%%%%%%%%%%%%%%%%%%%%%%%%%%%%%%%%%%%%%%%%%%%%%%%%%%%%%%%%%%%%%%%%%%%%%%%%%%%%%%%%%%%%%%%%%%%%%%%%
\section{Tour-5}
 The Euler's constant is defined as $$\gamma=\lim\limits_{n\to\infty}\gamma_{n},$$ where $$\gamma_{n}=1+\frac{1}{2}+\frac{1}{3}+\ldots\frac{1}{n}-\ln n.$$
In this visual tour, we will show that the existence of the  Euler's constant $\gamma$ and $\gamma\in(\frac{1}{2},1)$. Since
$$1-\gamma_{n}=\ln n-\frac{1}{2}-\frac{1}{3}-\ldots\frac{1}{n},$$ thus $1-\gamma_{n}$ can be described as the shaded area in the following figure.
\begin{figure}[h]
\includegraphics[scale=.5]{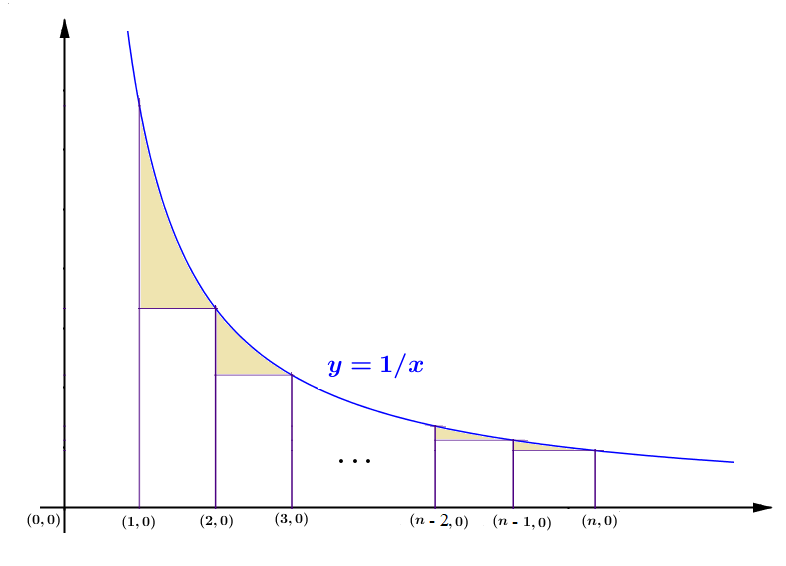}
\caption{$\ln n-\frac{1}{2}-\frac{1}{3}-\ldots-\frac{1}{n}=1-\gamma_{n}$}
\centering
\end{figure}
It is seen from the figure that $\{1-\gamma_{n}\}$ is strictly monotone increasing, and $1-\gamma_{n}>0.$ That is $\{\gamma_{n}\}$ is strictly monotone decreasing sequence and $\gamma_{n}$ is bounded above by $1$.
\newpage
Next, we define a sequence $\{A_{n}\}$, where
\begin{figure}[h]
\includegraphics[scale=.5]{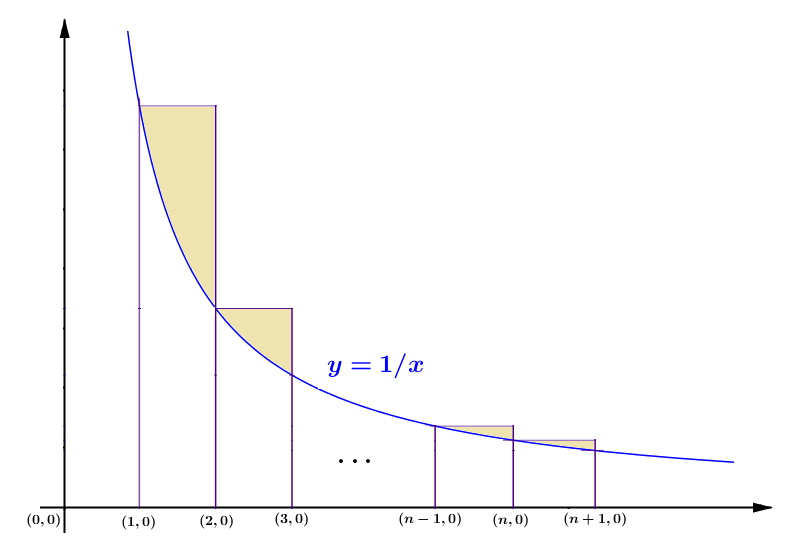}
\caption{$A_{n}=1+\frac{1}{2}+\frac{1}{3}+\ldots+\frac{1}{n}-\ln (n+1)$}
\centering
\end{figure}
$A_{n}$ is described by the shaded area in the following figure. Then it is seen that $\{A_{n}\}$ is strictly monotone increasing and $A_{n}>0.$
Since $$A_{n}=\gamma_{n}-\ln(n+1)+\ln n,$$ thus $$\gamma_{n}>\ln(n+1)-\ln n,$$ which means, by Figure 13, that
\begin{figure}[h]
\includegraphics[scale=.5]{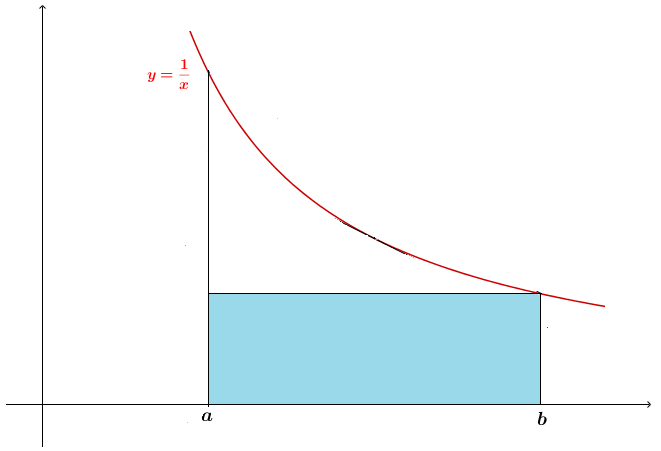}
\caption{$\ln b-\ln a >\frac{b-a}{b}$, \text{where} $b>a>0$.}
\centering
\end{figure}
$$\gamma_{n}>\ln(n+1)-\ln n>\frac{1}{n}>0,$$
i.e., $\gamma_{n}$ is bounded below by $0$. Thus the Euler's constant
$$\gamma=\lim\limits_{n\to\infty} \gamma_{n},$$ exist,  and $\gamma\in[0,1)$.
\newpage
 Next, we assume that $$\Gamma_{n}=1+\frac{1}{2}+\frac{1}{3}+\ldots\frac{1}{n}-\ln (n+\frac{1}{2}).$$ Thus, $$\Gamma_{n+1}-\Gamma_{n}=\frac{1}{n+1}+\ln(n+\frac{1}{2})-\ln(n+\frac{3}{2}).$$
Thus by the figure 14,
\begin{figure}[h]
\includegraphics[scale=.4]{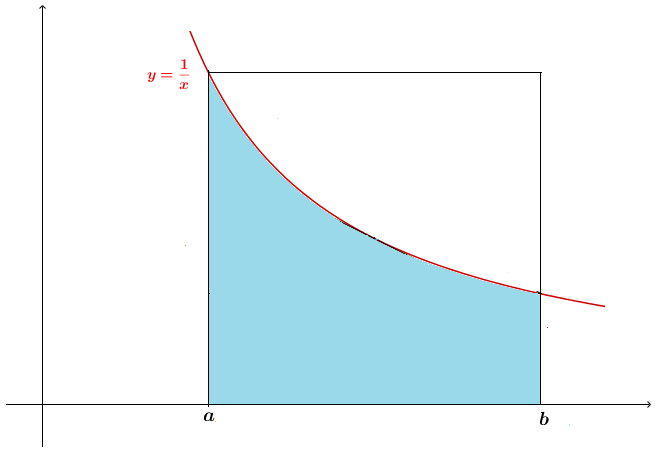}
\caption{$\ln b-\ln a <\frac{b-a}{a}$, \text{where} $b>a>0$.}
\centering
\end{figure}
$$\Gamma_{n+1}-\Gamma_{n}<\frac{1}{n+1}-\frac{1}{n+\frac{1}{2}}<0,$$
i.e., $\{\Gamma_{n}\}$ is strictly monotone decreasing sequence.\medbreak
%%%%%%%%%%%%%%%%%%%%%%%%%%%%%%%%%%%%%%%%%%%%%%%%%%%%%%%%%%%%%%%%%%%%%%%%%%%%%%%%%%%%%%%%%%%%%%%%%%%%%%%%%%%%%%%%%%%%%%%%%%%%%
Again, Figure 14 shows that
\beas \int_{n}^{n+\frac{1}{2}} \frac{1}{x}dx<\frac{1}{2n},\eeas
and, Figure 15 shows that
\begin{figure}[h]
\includegraphics[scale=.5]{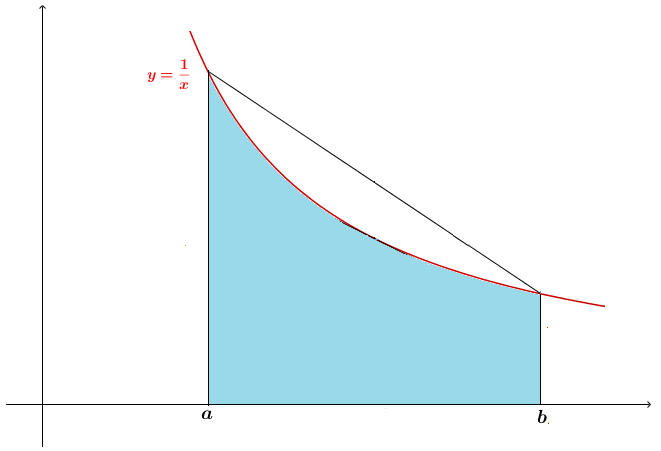}
\caption{$\ln b -\ln a< \frac{1}{2}(\frac{1}{a}+\frac{1}{b})\cdot(b-a).$}
\centering
\end{figure}
\beas \int_{1}^{n} \frac{1}{x}dx &=&\sum_{i=1}^{n-1} \int_{i}^{i+1} \frac{1}{x}dx\\
&<& \sum_{i=1}^{n-1} \frac{1}{2}(\frac{1}{i}+\frac{1}{i+1})\\
&=&1+\frac{1}{3}+\frac{1}{4}+\ldots+\frac{1}{n-1}+\frac{1}{2n}. \eeas
Thus \beas \ln\left(n+\frac{1}{2}\right)&=&\int_{1}^{n+\frac{1}{2}} \frac{1}{x}dx\\
&=& \int_{1}^{n} \frac{1}{x}dx+\int_{n}^{n+\frac{1}{2}} \frac{1}{x}dx\\
&<&1+\frac{1}{3}+\frac{1}{4}+\ldots+\frac{1}{n-1}+\frac{1}{n}, \eeas
Thus $\Gamma_{n}=1+\frac{1}{2}+\frac{1}{3}+\ldots\frac{1}{n}-\ln (n+\frac{1}{2})>\frac{1}{2}.$ Hence $\lim\limits_{n\to\infty} \Gamma_{n}$ exist and $\lim\limits_{n\to\infty} \Gamma_{n}\in [\frac{1}{2}, 1)$.\\
As  $\gamma_{n}-\Gamma_n=\ln (n+\frac{1}{2})-\ln n$, thus, applying Figures 13 and 14, we get
$$\frac{1}{2n+1}<\gamma_{n}-\Gamma_n<\frac{1}{2n},$$
i.e., $$\gamma=\lim\limits_{n\to\infty} \gamma_{n}=\lim\limits_{n\to\infty} \Gamma_{n}~~~\text{and}~~~~\gamma\in[\frac{1}{2},1).$$
%%%%%%%%%%%%%%%%%%%%%%%%%%%%%%%%%%%%%%%%%%%%%%%%%%%%%%%%%%%%%%%%%%%%%%%%%%%%%%%%%%%%%%%%
%%%%%%%%%%%%%%%%%%%%%%%%%%%%%%%%%%%%%%%%%%%%%%%%%%%%%%%%%%%%%%%%%%%%%%%%%%%%%%%%%%%%%%%%%%%%%%%%%%%%%%%%%%%%%%%%%%%%%%%%%%%%%%%%%%%%%%%%%%%%%%%%%%%%%%%%%%%%%%%

\end{document}